\documentclass[10pt]{amsart}
\usepackage{amsmath,amsfonts}

\newcommand{\sh}[1]{{#1}}

\newcommand{\C}{\mathbb{C}}
\newcommand{\Z}{\mathbb{Z}}
\newcommand{\CC}{\mathcal{C}}
\newcommand{\F}{\mathcal{F}}
\newcommand{\T}{\mathcal{T}}

\newcommand{\SL}{\mathop{\rm SL}\nolimits}
\newcommand{\Span}{\mathop{\rm Span}\nolimits}

\newcommand{\ord}{\mathop{\rm ord}\nolimits}

\newcommand{\up}{\mathop{\rm up}\nolimits}
\newcommand{\roof}{\mathop{\rm roof}\nolimits}

\newcommand{\ep}{\epsilon}
\newcommand{\Lam}{\Lambda}

\newcommand{\f}{{f}}
\newcommand{\e}{{e}}
\newcommand{\g}{\mathfrak{g}}
\newcommand{\n}{\mathfrak{n}}
\newcommand{\h}{\mathfrak{h}}
\renewcommand{\sl}{\mathfrak{sl}}
\newcommand{\gl}{\mathfrak{gl}}
\newcommand{\geqlex}{\stackrel{\rm lex}{\geq}}
\newcommand{\leqlex}{\stackrel{\rm lex}{\leq}}

\newcommand{\llex}{\stackrel{\rm lex}{<}}
\newcommand{\leqB}{\stackrel{\rm B}{\leq}}
\newcommand{\geqB}{\stackrel{\rm B}{\geq}}
 \newcommand{\hatE}{{\widehat{E}}}
\newcommand{\hatJ}{{\widehat{J}}}
\newcommand{\hatK}{{\widehat{K}}}

\newcommand{\hatk}{{\widehat{k}}}
\newcommand{\hatp}{{\widehat{p}}}
\newcommand{\hatq}{{\widehat{q}}}
\newcommand{\tp}{{\widetilde{p}}}
\newcommand{\Ct}{{\C[t^{\pm 1}]}}
\renewcommand{\mod}{\ \mathrm{mod}\ }
\newcommand{\height}{\mathop{\rm height}\nolimits}

\newcommand{\tJ}{{\widetilde J}}

\renewcommand{\Box}{\diamond}
\newcommand{\TT}{{
  \mbox{\scriptsize $\stackrel{\bf T}{}$ }}}
\newcommand{\bdot}{\mbox{\tiny $\bullet$}}
\newcommand{\wt}{\mathop{\rm wt}}

\newcommand{\lub}{\mathop{\rm lub}\nolimits} %V
\newcommand{\ceil}{\mathop{\rm ceil}\nolimits} %V
\newcommand{\bA}{\mathbf{a}} %V
\newcommand{\bB}{\mathbf{b}} %V
\newcommand{\bC}{\mathbf{c}} %V
\newcommand{\bD}{\mathbf{d}} %V
\newcommand{\bE}{\mathbf{e}} %V
\newcommand{\bF}{\mathbf{f}} %V
\newcommand{\bG}{\mathbf{g}} %V
\newcommand{\bH}{\mathbf{h}} %V
\begin{document}

\title{Standard Bases for Affine SL(n)-Modules}

\author{V. Kreiman}
\address{Department of Mathematics\\
Virginia Polytechnic Institute\\
Blacksburg, VA 24061}
\email{vkreiman@vt.edu}

\author[V. Lakshmibai]{V. Lakshmibai${}^{*}$}
\address{Department of Mathematics\\ Northeastern University\\
Boston, MA 02115}
\email{lakshmibai@neu.edu}
\thanks{${}^{*}$ Partially suported
by NSF grant DMS-0400679 and NSA-MDA904-03-1-0034.}

\author[P. Magyar]{P. Magyar${}^{\dag}$}
\address{Department of Mathematics\\
Michigan State University\\
East Lansing, MI 48824}
\email{magyar@math.msu.edu}
\thanks{${}^{\dag}$ Partially suported
by NSF grant DMS-0405948.}

\author[J. Weyman]{J. Weyman${}^{\ddag}$}
\address{Department of Mathematics\\
Northeastern University\\
Boston, MA 02115}
\email{j.weyman@neu.edu}
\thanks{${}^{\ddag}$ Partially
suported by NSF grant DMS-0300064.}

%\centerline{\Large\bf 
%Standard Bases for Affine SL(n)-Modules 
%}\vspace{1em}
%\centerline{ V.~Kreiman} 
%\centerline{ V.~Lakshmibai}
%\centerline{ P.~Magyar}
%\centerline{ J.~Weyman}
%\vspace{.5em}
%\centerline{Jan 1, 2005}
%\vspace{1em}

\begin{abstract}
We give an elementary and easily computable basis for
the Demazure modules in the basic representation of the affine Lie algebra $\widehat{\mathfrak{sl}}_n$ (and the loop group $\widehat\SL_n$).  A novel feature is that we define our basis ``bottom-up'' by raising each extremal weight vector, rather than ``top-down'' by lowering the highest weight vector. 

Our basis arises naturally from the combinatorics of its indexing set, which consists of certain subsets of the integers first specified by the Kyoto school in terms of crystal operators.  We give a new way of defining these special sets in terms of a recursive but very simple algorithm, the roof operator, which is analogous to the left-key construction of Lascoux-Schutzenberger.  The roof operator is in a sense orthogonal to the crystal operators.
\end{abstract}

\maketitle

\noindent
The most important representation of the affine Kac-Moody algebra
$\widehat\sl_n$\ (or of the loop group
$\widehat\SL_n$) is the basic representation $V(\Lam_0)$,
the highest-weight representation associated to
the extra node of the extended Dynkin diagram
$A_{n-1}^{(1)}$.  
The infinite-dimensional space
$V(\Lam_0)$ is filtered by the finite-dimensional 
Demazure modules $V_w(\Lam_0)$ for $w$ an element
of the affine Weyl group: these are modules for 
a Borel subgroup of the loop group.

There are several general constructions for irreducible representations and their Demazure modules, such as Lusztig's canonical basis \cite{lusztig}  
and Littelmann's contracting modules \cite{littelmann3}.
However, they are extremely difficult to compute explicitly, and even
the combinatorial indexing set for a basis is very intricate (see \cite{BZ}).  We will give an elementary and easily computable basis 
for $V(\Lam_0)$
and its Demazure modules.

We work inside the Fock space $\F$, an infinite wedge product which contains $V(\Lam_0)$, analogously to the space $\wedge^j\C^n$ which realizes a fundamental representation of $\SL_n\!\C$.
The Fock space has a natural basis indexed by certain infinite subsets of integers.
The combinatorial part of our problem amounts to defining which of these subsets will index basis elements of
$V_w(\Lam_0)$ for a given $w$.
We describe these special subsets in terms of a
recursive but very simple algorithm,
the roof operator on subsets.  This is analogous to the left-key construction of Lascoux-Schutzenberger \cite{LS}, which distinguishes the Young tableaux indexing a basis of a given Demazure module of $\SL_n\!\C$.  

The roof operator is more elementary (and much more efficient) than the crystal graph operators, and is in some sense orthogonal to them.  One may think of the roof operator as jumping across the crystal graph, moving each vertex down to an extremal weight vertex $w(\Lam_0)$, but {\it not} along edges of the crystal graph.

The combinatorics of the roof operator lead naturally to the definition of our standard basis, in analogy to the method of Raghavan-Sankaran \cite{raghavan}.  A novel feature is that we define our basis ``bottom-up'' by raising each extremal weight vector of $V(\Lam_0)$, rather than ``top-down'' by lowering the highest weight vector.  We prove linear independence of our basis by showing its triangular relationship to the natural basis of the Fock space.  We prove that our basis spans $V(\Lam_0)$ by showing that our special indexing subsets fill the crystal graph.

The paper is organized as follows.  In Section 1, we fix notation, define the roof operator and the standard basis,  state our main results, and point out related work.  In Section 2, we recall the basics of crystal graphs.  In Section 3, we prove the combinatorial comparison between the subsets distinguished by our roof operator and those in the crystal graph.  In Section 4, we prove the triangularity between bases in the Fock space.

\section{Main Results}

\noindent
Consider the complex untwisted affine Lie algebra of type $A_{n-1}^{(1)}$ :
$$\g=\widehat{\mathfrak{sl}}_n=
\sl_n(\Ct)\ \oplus\ \C K\ \oplus\ \C d
,$$ 
where $\sl_n(\Ct)$ denotes the traceless $n\times n$ matrices with entries in the Laurent polynomials $\Ct=\C[t,t^{-1}]$,\ $K$ is a central element of $\g$, and $d=t\frac{d}{dt}$ is a derivation (see \cite[Ch 7]{kac1}).  We have the Cartan decomposition $\g=\n\oplus\h\oplus\n_{-}$, where $\h$ is the Cartan subalgebra 
$$
\h=\bigoplus_{1\leq i\leq n-1}
\!\!\! \C(E_{ii}\sh-E_{i+1,i+1})\ 
\oplus\ \C K \ \oplus\ \C d\,;
$$ 
and $\n$ is the maximal nilpotent subalgebra
$$
\n:=\mathop{\bigoplus_{1\leq i<j\leq n}}_{k\geq 0}\!\!\!\C t^kE_{ij}\ \oplus\mathop{\bigoplus_{1\leq i<j\leq n}}_{k\geq 1}\!\!\!\C t^k E_{ji}
\ \oplus\
\mathop{\bigoplus_{1\leq i\leq n-1}}_{k\geq 1}t^k
(E_{ii}-E_{i+1,i+1})\,.
$$
Here $E_{ij}\in\gl_n(\C)$ denotes a coordinate matrix.
%We also have the quotient Lie algebra $\g'=\sl_n(\C)
%\otimes_\C \C[t,t^{-1}] = \sl_n(\C[t,t^{-1}])$.

Let $\Lam_0,\Lam_1,\ldots,\Lam_{n-1}$ be the fundamental weights of $\g$, and let $V(\Lam_m)$ be the level 1 irreducible $\g$-module with highest weight $\Lam_m$. (Thus, $V(\Lam_0)$ is the basic representation of $\g$.)\  Let us recall the construction of $V(\Lam_m)$ inside the fermionic Fock space $\F$ (cf.~\cite[Ch 14]{kac1}, \cite{kac2}).  Let $\C^\infty=\bigoplus_{i\in\Z} \C \ep_i$ be the $\C$-vector space with basis $\{\cdots,\ep_{-2},\ep_{-1},\ep_0,\ep_1,\ep_{2},\cdots\}$.
%, and let
%$$\begin{array}{rcl}
%\gl_\infty&=&\mathrm{End}_\C(\C^\infty)\\
%&=&
%\{\textstyle\sum_{i,j\in\Z}a_{ij} E_{ij}\ \mid \
%\text{for each $j$,\ $a_{ij}\neq 0$ for only finitely many %$i$}\}\,.
%\end{array}$$
%

Let $\T$ denote the collection of all subsets $J\subset\Z$ which are comparable to the non-positive integers $\Z_{\leq 0}$, meaning that $J\sh\setminus\Z_{\leq 0}$ and $\Z_{\leq 0}\sh\setminus J$ are both finite:
$$
\T=\{\,J\subset\Z\ \ \ \text{s.t.}\ \ \
|J\sh\setminus\Z_{\leq 0}|,\,|\Z_{\leq 0}\sh\setminus J|<\infty\,\}\,.
$$
We write such a set as:
$$
J=\{\cdots<j_{-2}<j_{-1}<j_0\}\,.
$$
Define the Fock space as the semi-infinite wedge product of $\C^\infty$:
$$\F=\wedge^{\!\infty/2}\,\C^\infty:=\bigoplus_{J\in\T}\C \ep_J\,,$$ 
the $\C$-vector space with basis elements:
$$
\ep_J:=\,\cdots\sh\wedge\ep_{j_{-2}}\sh\wedge\ep_{j_{-1}}\sh\wedge\ep_{j_{0}}\,.
$$  
Thus, the $J\in\T$ play the role of tableaux indexing the basis vectors of the Fock space.

For $i,j\in\Z$, let $E_{ij}'$ denote a coordinate 
matrix acting on $\C^\infty$ by $E_{ij}'(\ep_j)=\ep_i$,\ \, $E_{ij}'(\ep_k)=0$ for $k\neq j$.
Then for $i<j$,\ $E_{ij}'$ acts on the Fock space in the expected way:
$$  
E_{ij}'(\ep_J) = \left\{\begin{array}{cl}
\pm \ep_{J\setminus j\cup i}&\text{if } j\in J,\ i\not\in J\\[.5em]
0& \text{otherwise,} 
\end{array}\right.$$
Here we denote:
$$
J\setminus j\cup i\,:=\,(\,J\sh\setminus\{j\}\,)\,\cup\,\{i\}\,,
$$
the operation which moves the element $j\in J$ to the vacant position $i\not\in J$;
and $\pm=(-1)^\ell$ with $\ell={|J\cap[i,j]|-1}$, the sign of the permutation needed to sort the wedge factors of $\ep_{J\setminus j\cup i}$ into increasing order.

We let: 
$$
\hatE_{pq}:=\sum_{k\in\Z} E_{p+nk,q+nk}'\,,
$$
which is a well-defined operator on $\F$.
Now, if $i<j$ or $k>0$, we let $t^k E_{ij}$ act on $\F$ by the operator
$\hatE_{pq}$, where $p=i-nk$, $q=j$:
$$
t^k E_{ij}=\hatE_{i-nk,\,j}:\F\to\F\,.
$$
This defines the action\begin{footnote}{
This action arises naturally 
if we identify the free $\Ct$-module
$\Ct^n=\bigoplus_{i=1}^n\Ct \ep_i$ with the $\C$-vector space
$\C^\infty=\bigoplus_{j\in\Z}\C\ep_j$
via: $t^k \ep_i\leftrightarrow\ep_{i-nk}\,.$
This gives an embedding $\gl_n(\Ct)\subset\gl(\C^\infty)$,
so that the natural action of the upper triangular part of $\gl(\C^\infty)$ on the Fock space restricts to the specified action of $\n\subset\gl_n(\Ct)$.  
However, this gives only a projective representation of the entire $\gl_n(\Ct)$, which then lifts to a true representation of the central extension $\widehat{\mathfrak{gl}}_n$.}
\end{footnote}of $\n$ on $\F$, and we can similarly define the action of $\n_{-}$ and $\h$.
Indeed, the Chevalley generators of $\n_{-}$ are:
$F_i=E_{i+1,i}=\hatE_{i+1,i}$ for $i=1,\ldots,n\sh-1$ and $F_0=t^{-1}E_{1,n}=\hatE_{n+1,n}$.

Now let $L_m:=\Z_{\leq m}\in \T$.
It is well known that the $U(\g)$-span of
the highest-weight vector $\ep_{L_m}$ is an irreducible $\g$-module:
$$
U(\g)\cdot \ep_{L_m} = U(\n_{-})\cdot \ep_{L_m} \cong
V(\Lam_m)\,,
$$
where we define $\Lam_{m}:=\Lam_{(m \mod n)}$.

Recall that we can realize the Weyl group $W$ of $\g$ as a permutation group on $\Z$.  Indeed, we can write the simple reflection $s_i:\Z\to\Z$ as a product of commuting transpositions:
%$$
%s_i(j)=\left\{\begin{array}{cl}
%j+1 & \text{if}\ \ j\equiv i\mod n\\
%j-1 & \text{if}\ \ j\equiv i\sh+1\mod n\\
%j & \text{otherwise},
%\end{array}\right.$$
$$
s_i:=\prod_{k\in\Z} (i\sh+nk,\,i\sh+1\sh+nk)\ ,
$$
so that $s_i(i')=i'\sh+1$ whenever $i'\equiv i\mod n$. 
Then $W=\langle s_0,\ldots,s_{n-1}\rangle$ is the corresponding Coxeter group.

The Weyl group $W$ acts on $\T$ via 
$w(J):=\{w(j)\}_{j\in J}$. 
Indeed, the extremal weight vectors of $V(\Lam_m)\subset\F$ are just $\ep_J$ for $J=w(L_m)$.  
Equivalently, a basis vector $\ep_J$ is an extremal weight vector whenever $J$ is {\it $n$-stable}:  that is, whenever $j-n\in J$ for all $j\in J$.  
We define the {\it parabolic Bruhat order} between $K=\{\cdots\sh<k_{-1}\sh<k_0\}$ and $J=\{\cdots\sh<j_{-1}\sh<j_{0}\}$ as:
$$
K\leqB J \quad\Longleftrightarrow\quad
\left\{\begin{array}{cl}
k_i\leq j_i&\text{ for all } i\\
k_i=j_i&\text{ for all } i\sh\ll 0
\end{array}\right.
$$ 
This induces an order on the $n$-stable $J=w(L_m)$ which is consistent with the usual Bruhat order on $w\in W$.

The {\it Demazure modules} \cite{demazure} of $V(\Lam)$ are the $\n$-modules obtained by raising the extremal weights:
$$
V_w(\Lam_m)\ \cong\ U(\n)\cdot \ep_{w(L_m)}\,.
$$
We can get the same modules also by lowering the highest weight:
$$
V_w(\Lam_m)=\Span_\C\{F_{i_1}^{k_1}\cdots F_{i_t}^{k_t}
\ep_{L_m}\
\mid k_1,\ldots,k_r\geq 0\}\,,
$$
where $w=s_{i_1}\cdots s_{i_t}$ is a reduced word.

Next we describe the sets $J\in\T$ which index basis vectors of $$V_w(\Lam_m)\subset V(\Lam_m)\subset\F\,.$$
Let us say that a set $J\in\T$ is 
{\it n-bounded} if $j_i-j_{i-1}\leq n$ for all $i$. 
Also, we define the {\it order} of a set $J$ by:
$\ord(J):=|J\sh\setminus \Z_{\leq 0}|-|\Z_{\leq 0}\sh\setminus J|$ ; equivalently, $\ord(J)=m$ means that
$j_i=m+i$ for all sufficiently large negative $i$.
Now let
$$\begin{array}{rcl}
\CC(L_m)&:=&
\{J\in \T\ \mid\ \ord(J)=m\,\text{ and $J$ is $n$-bounded}\}
\\[.3em]
&=&\left\{J=\{\cdots\sh<j_{-2}\sh<j_{-1}\sh<j_{0}\}\subset\Z\,\left|\begin{array}{c}
j_i=m+i\text{ for }i\ll 0\\
j_i-j_{i-1}\leq n\text{ for all }i
\end{array}\right.\right\}
\end{array}$$
(The reader should be aware of a frequently used alternative notation in terms of ``colored Young diagrams'' instead of subsets.\!\!
\footnote{In \cite{FLOTW} and related literature, the basis
 of $V(\Lam_m)$ is indexed by the set of all
partitions $\lambda=(\lambda_1\sh\geq\lambda_2\sh\geq\cdots)$ with $\lambda_i\geq 0$, \ $\lambda_i=0$ for $i\gg0$, and 
$\lambda_{i+1}\sh-\lambda_i\leq n\sh-1$.  Namely, a set $J=\{\cdots\sh< j_{-1}\sh<j_0\}$ of order $m$ corresponds to $\lambda$ with $\lambda_{i+1}=m-i-j_{-i}$.  It is useful to picture $\lambda$ as a Young diagram colored with a mod-$n$ checkerboard pattern: square $(i,j)$ has color $i-j\in\Z/n\Z$. 
})

We can give $\CC(L_m)$ a crystal graph structure by defining the crystal lowering operators $\f_i$
for $i=0,\ldots,n\sh-1$, as recalled in Section 2 below. 
If it is defined, the crystal operator $f_i$ on a set $J$ picks out a certain element $r\in J$ with $r\equiv i\mod n$, and replaces it with $r\sh+1\equiv i\sh+1\mod n$: that is,
$f_i(J)=J\setminus r\cup(r\sh+1)\,.$ 
%The raising operator $e_i$ picks out a certain element $q\in J$ with $q\equiv i\sh+1\mod n$, and replaces it by $q\sh-1\equiv i\mod n$.
%We have $K=f_i(J)$ if and only if $e_i(K)=J$.
%
We define the Demazure crystal as:
$$\CC_w(L_m)=\{ \f_{i_1}^{k_1}\cdots \f_{i_t}^{k_t} L_m\mid k_1,\ldots,k_t\geq 0\},$$
where $w=s_{i_1}\cdots s_{i_t}$ is again a reduced word.

Our first theorem is a simpler description of the sets $J$ in this Demazure crystal, in analogy with the ``left key'' algorithm of Lascoux-Schutzenberger \cite{LS}.
If $J$ is $n$-bounded but not $n$-stable, define the following up-operation (which is different from the crystal operators): 
$$\begin{array}{c}
\up(J):=J\setminus p\cup q\quad \text{where:}\\[.5em]
p:=\max\{p'\ \mid\  p'\in J,\ \,p'\sh-n\not\in J\},\quad
\\[.5em]
q:=\min\{q'> p\ \mid\ q '\not\in J,\ \ q'\sh-n\in J,\ \
q'\not\equiv p \mod n\}\,.
\end{array}$$
To rephrase this in words, define a {\it seam} as a maximal arithmetic progression $S=\{\cdots<j\sh-2n<j\sh-n<j\}$ contained in $J$.  We call the vacant position $j\sh+n\not\in J$ the {\it tight end} of $S$\,; and if $S$ is finite, we call the minimal element $p\in S$ the {\it loose end} of $S$.  The up-operation moves $p\in J$ to $q\not\in J$, where $p$ is the maximal loose end in $J$, part of a seam $S=\{p,p\sh+n,\ldots\}$, and $q>p$ is the tight end of a different seam, the minimal such tight end.  See the examples below.

Iterating the up-operation ``pulls out" this seam, distributing all the elements of $S$ to the tight ends of different seams; and then the operation starts on another seam.  After each seam is pulled out, the number of loose ends decreases by one.  Once all the finite seams of $J$ are pulled out, the result is an $n$-stable set which we call the {\it roof} of $J$:
$$
\roof(J):=\up^\ell(J)=\up(\cdots\up(J)\cdots)=y(L_m)
\quad\text{for some}\ y\in W\,.
$$
{\bf Theorem 1}\quad {\it
Let $\CC_w(L_m)$ be the Demazure crystal generated from the highest weight $L_m$ according to a reduced word for $w\in W$.  Then:
$$\begin{array}{rcl}
\CC_w(L_m) &=&\{J\in\CC(L_m) \ |\ \roof(J)\leq w(L_m)\}
%\\[.5em]
%&=&\{J\subset\Z \ |\  \ord(J)=m,\ 
%J\ \text{\rm is $n$-bounded},\ \roof(J)\leq w\,\}
\end{array}$$
}
\\[-1.2em]
This gives a highly efficient algorithm for testing the membership of $J$ in $\CC_w(L_m)$.  (Later in this section we give a corresponding algorithm for generating all $J\in\CC_w(L_m)$.)

Next we give elementary bases of $V(\Lam_m)$ and its dual which are compatible with the Demazure modules, in analogy to the construction of Raghavan-Sankaran \cite{raghavan} (generalized by Littelmann \cite{littelmann3}).
\\[1em]
{\bf Theorem 2}\quad {\it
(i) Given $J\in\CC_w(L_m)$, suppose $\up^i(J)=\up^{i-1}(J)\setminus p_i\cup q_i$ for $i=1,\ldots,\ell$, and $\roof(J)=\up^\ell(J)=y(L_m)$.
%Here $p_i<q_i$.
Define $$ v_J:= 
\hatE_{p_1,q_1}\cdots\hatE_{p_\ell,q_\ell}\ep_{y(L_m)}\,.$$
Then the irreducible highest-weight module $V(\Lam_m)$,
a submodule of the Fock space $\F$, has basis $\{v_J\mid J\in\CC(L_m)\}$; and the Demazure module $V_w(\Lam_m)$ has basis $\{v_J\mid J\in\CC_w(L_m)\}$.  
\\[.5em]
(ii) The irreducible lowest-weight module $V(\Lam_m)^*$, 
a quotient of the dual Fock space $\F^*$, has basis
$\{\ep_J^*\mid J\in\CC(L_m)\}$; and the dual Demazure module $V_w(\Lam_m)^*$ has basis
$\{\ep_J^*\mid J\in\CC_w(L_m)\}$.
Here $\ep_J^*$ denotes a dual basis vector of $\F^*$ restricted to $V(\Lam_m)$ or to $V_w(\Lam_m)$ respectively.
\\[.5em]
}
In geometric terms, the functions $\ep_J^*$ can be considered as Plucker coordinates on the affine Grassmannian embedded in the infinite projective space $\mathbb{P}(\F)$.

The vectors $v_J$ possess a triangularity property
with respect to the standard basis of the Fock space.
Define {\it lexicographic order} on sets 
$K,J$ as follows:
$$
K\llex J\quad\Longleftrightarrow\quad\left\{\begin{array}{cl}
k_N<j_N&\text{for some $N$}\\
k_i=j_i &\text{for all $i<N$\,.}
\end{array}\right.
$$
{\bf Proposition 3}\quad {\it
Let us write: 
$$v_J=\sum_K a_K^J \ep_K\quad \text{with coefficients}\quad\ a_K^J\in\C\,.$$
\\[-1em]
(i) We have a non-zero coefficient 
$a_K^J\neq 0$ only if $K\leqlex J$.\\
(ii) We have a non-zero leading coefficient $a_J^J\neq 0$  for every $J$, given explicitly as follows.
If $J$ is $n$-stable, then $a_J^J=1$.  
If $J$ is $n$-bounded but not $n$-stable, suppose that $\up^i(J)=\up^{i-1}(J)\setminus p_i\cup q_i$, and $t$ is maximal such that
$p_1\equiv\cdots\equiv p_t\mod n\,.$ 
Define $\mu_d:=\#\{\,i\leq t\mid q_i-p_i=d\,\}\,$ and  $\tJ:=\up^t(J)$. Then:
$$
a_J^J=\pm\left(\prod_{d\geq 1}\mu_d!\right) a_\tJ^\tJ\,.
$$
}
\\[0em]
In part (ii), note that the sequence $S=\{p_1<\ldots<p_t\}$ is actually the first seam of $J$ pulled out by the roof algorithm:\ 
$S=\{p\,,\,p\sh+n,\ldots, p\sh+n(t\sh-1)\}$.
Iterating part (ii), we get a combinatorial formula for
the leading coefficient $a_J^J$ of each $v_J$ depending only on the sequences $p_1,\ldots,p_\ell$ and $q_1,\ldots,q_\ell$ in the roof algorithm.  
\\[1em]
{\bf Example}\quad Let $n=5$ and let:
$$
J:=\{\ldots,-4,-3,-2,-1, 0,3,4,7,10,12,14,17,18,23,27,32,33,35,37\}\,.
$$
Then $J\in\CC(L_m)$ for $m=14$, since
$L_0\subset J$ and $|J\sh\setminus L_0|=14$, so that
$\ord(J)=\ord(L_0)\sh+14=14$.
We sort $J$ into its residue classes mod $n$ to show the seam structure. We mark the maximal loose end with boldface, and the tight end used by the up-operation with $\TT$. 
$$\begin{array}{r@{\,}c@{\,}l}
J&=&\!\!
\mbox{\footnotesize$\left[\begin{array}{@{\ }c@{\ \,}c@{\ \,}c@{\ \,}c@{\ \,}c@{\ \,}c@{\ \,}c@{\ \,}c@{\ \,}c@{\ \,}c@{\ \,}c@{\ \,}c@{\,}}
\cdots&-4&\cdot&\cdot&\cdot&\cdot&\cdot&\cdot&\cdot&\cdot\\
\cdots&-3&\cdot&7&12&17&\cdot&27&32&37\\
\cdots&-2&3&\cdot&\cdot&18&23&\cdot&33&\cdot\\
\cdots&-1&4&\cdot&14&\cdot&\cdot&\cdot&\cdot&\cdot\\
\cdots&\ \ 0&\cdot&10&\cdot&\cdot&\cdot&\cdot&{35}&\cdot&
\end{array}\!\!\!\right]$}
=
\mbox{\footnotesize$\left[\begin{array}{@{\ }c@{\ \,}c@{\ \,}c@{\ \,}c@{\ \,}c@{\ \,}c@{\ \,}c@{\ \,}c@{\ \,}c@{\ \,}c@{\ \,}c@{\ \,}c@{\,}}
\cdot&\cdot&\cdot&\cdot&\cdot&\cdot&\cdot&\cdot\\
\cdot&7&12&17&\cdot&27&32&37\\
3&\cdot&\cdot&18&23&\cdot&33&\TT\\
4&\cdot&14&\cdot&\cdot&\cdot&\cdot&\cdot\\
\cdot&10&\cdot&\cdot&\cdot&\cdot&\mathbf{35}&\cdot&
\end{array}\!\!\!\!\right]$}
\\[2.5em]
&\stackrel{\up}{\to}&
\!\!\mbox{\footnotesize$\left[\begin{array}{@{\ }c@{\ \,}c@{\ \,}c@{\ \,}c@{\ \,}c@{\ \,}c@{\ \,}c@{\ \,}c@{\ \,}c@{\ \,}c@{\ \,}c@{\ \,}c@{\,}}
\cdot&\cdot&\cdot&\cdot&\cdot&\cdot&\cdot&\cdot&\cdot\\
\cdot&7&12&17&\cdot&27&32&37&\TT\\
3&\cdot&\cdot&18&23&\cdot&\mathbf{33}&38&\cdot\\
4&\cdot&14&\cdot&\cdot&\cdot&\cdot&\cdot&\cdot\\
\cdot&10&\cdot&\cdot&\cdot&\cdot&\cdot&\cdot&\cdot
\end{array}\!\!\!\right]$}
\stackrel{\ \up^{\mbox{\scriptsize 2}}}{\longrightarrow}
\mbox{\footnotesize$\left[\begin{array}
{@{\ }c@{\ \,}c@{\ \,}c@{\ \,}c@{\ \,}c@{\ \,}c@{\ \,}c@{\ \,}c@{\ \,}c@{\ \,}c@{\ \,}c@{\ \,}c@{\ \,}c@{\,}}
\cdot&\cdot&\cdot&\cdot&\cdot&\cdot&\cdot&\cdot&\cdot&\cdot\\
\cdot&7&12&17&\cdot&\mathbf{27}&{32}&{37}&{42}&{47}\\
3&\cdot&\cdot&18&23&\TT&\cdot&\cdot&\cdot&\cdot\\
4&\cdot&14&\cdot&\cdot&\cdot&\cdot&\cdot&\cdot&\cdot\\
\cdot&10&\cdot&\cdot&\cdot&\cdot&\cdot&\cdot&\cdot&\cdot
\end{array}\!\!\!\right]$}
\\[2.5em]
&\stackrel{\ \up^{\mbox{\scriptsize 5}}}{\longrightarrow}&
\!\!\mbox{\footnotesize$\left[\begin{array}
{@{\ }c@{\ \,}c@{\ \,}c@{\ \,}c@{\ \,}c@{\ \,}c@{\ \,}c@{\ \,}c@{\ \,}c@{\ \,}c@{\ \,}c@{\ \,}c@{\ \,}c@{\,}}
\cdot&\cdot&\cdot&\cdot&\cdot&\cdot&\cdot&\cdot&\cdot&\cdot\\
\cdot&7&12&17&\cdot&\cdot&\cdot&\cdot&\cdot&\\
3&\cdot&\cdot&\mathbf{18}&23&28&33&38&43&48\\
4&\cdot&14&\TT&\cdot&\cdot&\cdot&\cdot&\cdot&\cdot\\
\cdot&10&\cdot&\cdot&\cdot&\cdot&\cdot&\cdot&\cdot&\cdot
\end{array}\!\!\!\right]$}
\!\!\stackrel{\ \up^{\mbox{\scriptsize 16}}}{\longrightarrow}\!\!\!\!
\mbox{
\footnotesize$\left[\begin{array}
{@{\ }c@{\ \,}c@{\ \,}c@{\ \,}c@{\ \,}c@{\ \,}c@{\ \,}c@{\ \,}c@{\ \,}c@{\ \,}c@{\ \,}c@{\ \,}c@{\ \,}c@{\ \,}c@{\,}}
\cdot&\cdot&\cdot&\cdot&\cdot&\cdot&\cdot&\cdot&\cdot&\cdot&\cdot&\cdot&\cdot\\
\cdot&\!\!\mathbf{7}&\!\!12&17&22&27&32&37&42&47&52&57&62\\
3&\TT&\cdot&\cdot&\cdot&\cdot&\cdot&\cdot&\cdot&\cdot&\cdot&\cdot\\
4&\cdot&\cdot&\cdot&\cdot&\cdot&\cdot&\cdot&\cdot&\cdot&\cdot&\cdot&\cdot\\
\cdot&\cdot&\cdot&\cdot&\cdot&\cdot&\cdot&\cdot&\cdot&\cdot&\cdot&\cdot&\cdot
\end{array}\!\!\!\right]$}
\\[2.5em]
&\stackrel{\ \up^{\mbox{\scriptsize 12}}}{\longrightarrow}&
\!\!\mbox{
\footnotesize$\left[\begin{array}
{@{\ }c@{\ \,}c@{\ \,}c@{\ \,}c@{\ \,}c@{\ \,}c@{\ \,}c@{\ \,}c@{\ \,}c@{\ \,}c@{\ \,}c@{\ \,}c@{\ \,}c@{\ \,}c@{\,}}
\cdot&\cdot&\cdot&\cdot&\cdot&\cdot&\cdot&\cdot&\cdot&\cdot&\cdot\\
\cdot&\cdot&\cdot&\cdot&\cdot&\cdot&\cdot&\cdot&\cdot&\cdot&\cdot&\cdot&\cdot\\
3&8&13&18&23&28&33&38&43&48&53&58&63\\
4&\cdot&\cdot&\cdot&\cdot&\cdot&\cdot&\cdot&\cdot&\cdot&\cdot\\
\cdot&\cdot&\cdot&\cdot&\cdot&\cdot&\cdot&\cdot&\cdot&\cdot&\cdot
\end{array}\!\!\!\right]$}
=\roof(J)=y(L_{14})\,.
\end{array}$$
We thus have $(p_1,q_1)\sh=(35,38)$,\ 
$(p_2,q_2)\sh=(33,42)$,\ 
$(p_3,q_3)\sh=(38,47)$,\ldots, and:
$$\begin{array}{rcl}
v_J&=& \hatE_{35,38}\,\hatE_{33,42}\,\hatE_{38,47}
\cdots \ep_{y(L_{14})}\\[.5em]
&=&\hatE_{5,8}
\hatE_{3,12}^{\,2}\,
\hatE_{2,3}^{\,5}\,
\hatE_{3,4}^{\,7}\,
\hatE_{4,7}^{\,8}\,
\hatE_{5,57}\,
\hatE_{2,3}^{\,12}
\ \ep_{y(L_14)}\,.
\end{array}$$
Since $J$ has seven loose ends, we must apply Proposition 3(ii) seven times to compute:
$
a_J^J=\pm\ 1!\cdot 2!\cdot  5!\cdot  7!\cdot  8!\cdot
1!\cdot  1!\cdot  12!\,.
$
(Here we have only one factorial for each seam, though in general there will be several.)

To determine a reduced decomposition for the Weyl group element $y$, we start with the extremal weight $K=y(L_{14})$ and perform the simple reflection: 
$K\mapsto s_{r}(K)$,
where 
$$
r:=\min\,\{\,k\!\not\in\! K\ \mid\ k\sh+1\!\in\! K\,\}\,
$$
is the minimal ``hole'' of $K$, and $s_r:=s_{(r\mod n)}$.  This will always give
$K\stackrel{\mathrm{B}}{>}s_r(K)$, and iterating the operation produces a canonical reduced word for $y$.
Indeed, $$
y(L_{14})=s_{2}\,s_{1}\,s_{3}\,s_{2}\,s_{0}\,
(s_{4}s_{3}\,s_{2}\,s_{1}\,s_{0})^{11}\,s_{4}\, L_{14}
\ .$$
See also the Example in the next section.
\quad$\Box$
\\[1em]
{\bf Example}\quad Let $n=2$.  Then for fixed $m$, the Bruhat order on the sets $w(L_m)$ reduces to a linear order: for example, 
$$L_0 = s_1(L_0)
\stackrel{\rm B}{<}
s_{0}(L_0)
\stackrel{\rm B}{<}
s_{1}s_{0}(L_0)
\stackrel{\rm B}{<}
s_{0}s_{1}s_{0}(L_0)
\stackrel{\rm B}{<}\cdots\,.$$
For $J$ an $n$-bounded set with order $m$, the roof operation reduces to:
$$\begin{array}{rcl}
\roof(J)&=&\min\{\,w(L_m)\,\mid\, w\in W,\ J\leqB w(L_m)\,\}\\[.5em]
&=&L_a\cup\{a\sh+2,a\sh+4,\ldots,a\sh+2k\}\,,
\end{array}$$
where $a=\max\{a'\mid L_{a'}\subset J\}$ and $k=|J\setminus L_a|$.  That is, the Demazure crystal is simply $\CC_w(L_m)=\{J\in\CC(L_m)\,\mid\,J\leqB w(L_m)\}$.
This case is further considered in the context of completely integrable lattice models in \cite{FMO}.  \quad$\Box$
\\[1em]
{\bf Example}\quad Generalizing the previous case, let $n$ be arbitrary and suppose $w(L_m)$ is of the form:
$$
w(L_m) = L_a\cup
\{a\sh+n, a\sh+2n, \ldots, a\sh+kn\}
$$ 
for some $a$ and $k = m-a$.  Then for any $J\leqB w(L_m)$, we have $\roof(J)\leqB w(L_m)$, so that again $\CC_w(L_m)=\{J\in\CC(L_m)\,\mid\,J\leqB w(L_m)\}$.  This case, which is considered in \cite{MW}, is exceptional: in general, it often happens that $J\leqB w(L_m)$, but $\roof(J)\not\leqB w(L_m)$.
\quad$\Box$
\\[1em]
Next we consider the modifications which must be made to our theory to generalize it to positive characteristic.
Since the leading coefficients $a_J^J$ are not necessarily $\pm 1$, the vectors $v_J$ could become linearly dependent if we work over $\Z$ and then reduce modulo a prime. To define a characteristic-free basis $\{\,v^{\,\prime}_J\mid J\in\CC_w(L_m)\,\}$, we start with $v^{\,\prime}_J:=\ep_J$ for $n$-stable $J$, then for a general $n$-bounded $J$ we define:
$$
v^{\,\prime}_J\ \ :=\ \ \frac
{\hatE_{p_1,q_1}\cdots\hatE_{p_t,q_t}}
{\prod_{d\in\Z}\mu_d!}\cdot v^{\,\prime}_\tJ
\ \ =\ \
\mathop{\prod_{d\in\Z}}_{d\not\equiv 0}\frac{\hatE_{p,p+d}^{\,\mu_d}}{\mu_d!}\cdot v^{\,\prime}_\tJ\ ,
$$
where $t$ is maximal such that $p:=p_1\equiv\cdots\equiv p_t\mod n$.  The second equality follows because $d_i:=q_i-p_i\not\equiv 0\mod n$,  so all the operators $\hatE_{p,p+d}$ commute with each other.  The basis $\{v^{\,\prime}_J\}$ clearly lies in the Kostant $\Z$-form  of the Demazure module $V_w(\Lam_m)$, and it has leading coefficients $\pm1$, so it reduces to a basis over an arbitrary field.  (Cf. \cite[Ch. 26]{humphreys}.)

Theorem 1 also gives an alternative ``bottom-up'' algorithm to generate $\CC_w(L_m)$, as opposed to the ``top-down'' definition in terms of crystal lowering operators.  We write: 
$$\begin{array}{rcl}
\CC_{=y}(L_m)&:=&\{J\in\CC(L_m)\,\mid\,\roof(J)=y\}\\[.5em]
&=&\{y(L_m)\}\,\cup\,\up^{-1}(y(L_m))\,\cup\,
\up^{-1}\up^{-1}(y(L_m))\,\cup\, \cdots
\end{array}$$
where $\up^{-1}(\hatJ)$ means the set of all $J$ such that $\up(J)=\hatJ$. To compute this for any given $\hatJ\in\CC(L_m)$, we first find $\tp<\hatp$\,, the two maximal loose ends of $\hatJ$ (with one or both possibly $=-\infty$).  Next we choose any $q>\hatp\sh-n$ such that $q\sh+n$ is the tight end of a seam $S\subset\hatJ$ of length $|S|\geq 2$,
% $q\in\hatJ,\ q\sh+n\not\in\hatJ$ 
% (a maximal element of a seam of $\hatJ$),
and we let $\hatq$ be the maximal tight end of $\hatJ$ 
less than $q$.  Finally, we define:
$$
P(q):=\left\{\, p\ \left|\begin{array}{c}
p\sh-n,p\not\in\hatJ\\[.3em]
\max(\hatp,\hatq\,)<p<q
\end{array}\right.\right\}
\,\cup\,
\left\{\, p=\hatp\sh-n\ \left|\begin{array}{c}
p\sh-n\not\in\hatJ\\[.3em]
\max(\tp,\hatq\,)<p
\end{array}\right.\right\}\,.
$$
Then we have: 
$$
\up^{-1}(\hatJ)=
\left\{J:=\hatJ\setminus q\cup p\ \left|\
\begin{array}{c}
q\sh-n,q\in\hatJ,\,\ q\sh+n\not\in\hatJ\\[.3em]
q>\hatp\sh-n,
\,\ p\in P(q)
\end{array}\right.\right\}\,.
$$
Applying this to all $y\leq w$, we generate all $J\in\CC_w(L_m)$.

We will prove Theorem 1 in Section 3 and Proposition 3 in Section 4.  Theorem 2 is a corollary of these, as follows.
By Theorem 1 and the definitions, we have:
$$
V':=\Span_\C\{\,v_J\mid \roof(J)\leq w\,\}
=\Span_\C\{\,v_J\mid J\in\CC_w(L_m)\,\}
\subset V_w(\Lam_m)\,.
$$
Proposition 3 implies that the $v_J$ are linearly independent vectors in $\F$ (since they are triangular with respect to the standard  basis 
$\{\ep_J\}$ ), so that 
$\dim_\C V'= |\CC_w(L_m)|\,;$
but it is well known from crystal graph theory  (Section 2 below) that $|\CC_w(L_m)|=\dim_\C V_w(\Lam_m)$, so that
$V'=V_w(\Lam_m)$. This shows Theorem 2 for the Demazure module $V_w(\Lam_m)$, and the claims for the irreducible module and the dual modules follow trivially.

We comment on related work which is closest to our point of view.  The pioneering paper \cite{DJKMO} by Date, Jimbo, Kuniba, Miwa, and Okado of the Kyoto school defined the tableaux $\CC(L_m)$ for $V(\Lam_m)$ (\,in fact for all $V(\ell\Lam_m)$\,), and the crystal graph structure was first defined by Misra, Miwa, Jimbo, et al.~in \cite{MM}, \cite{jimbo}.  Certain Demazure crystals $\CC_w(L_m)$ were considered by Kuniba, Misra, Miwa, Uchiyama and others in \cite{KMOTU},\cite{KMOU},\cite{FMO},\cite{MW}. 
A useful survey of related work is \cite{FLOTW}, and \cite{kashiwara} is a fundamental reference.
\\[1em]
{\bf Notation}\quad  For a set $J\subset\Z$, we define:
$$
J^{\equiv i}:=\{j\in J\,\mid\, j\equiv i\mod n\}\,,
\qquad 
J_{<r}:=\{j\in J\,\mid\, j<r\}\,.
$$
Similarly for\ $J_{>r}$ , for\ \
$J^{\equiv i}_{>r}:=J^{\equiv i} \cap J_{>r}$ , 
for\ \ ${}_{q\leq}J_{\leq r}:=J_{\geq q}\cap J_{\leq r}$ ,
etc.

\section{Crystal Operators}

In this section, we review the necessary facts about the crystal raising and lowering operators acting on $\CC(L_m)$.\,\begin{footnote}{These operators are sometimes encoded in the {\it crystal graph} having vertices $J\in\CC(L_m)$ and $i$-colored edges $J\stackrel{i}{\to}f_iJ$.}
\end{footnote}These operators were first defined in our case by the Kyoto school \cite{jimbo}, and they can also be derived from Littelmann's path model (as modified for semi-infinite paths in \cite{magyar}).   The crystal operators are basically different from the up-operation: indeed, by Theorem 1 the two are in some sense transversal to each other.

If it is defined, the lowering operator $f_i$ for $i=0,1,\ldots n\sh-1$ acting on a set $J\in\T$ picks out a certain element $r\in J$ with $r\equiv i\mod n$, and replaces it with $r\sh+1\equiv i\sh+1\mod n$.
(We say that $f_i(J)$ is ``lower'' than $J$ because it is farther from the highest-weight element $L_m$.)\ \ 
Similarly, the raising operator $e_i(J)$ picks out a certain element $r'\in J$ with $r'\equiv i\sh+1$ and replaces it with $r'\sh-1\equiv i$.  We have:\ \ $f_i(J)=J' \,\Longleftrightarrow\, J=e_i(J')$.
\\[1em]
{\bf Definition}\ \ {\it Given $J\in\T$.\\
(i) Let $$
R:=\{\,r\ \text{\ \ s.t.\ \ for all $k\geq r$}\,,\ \
|{}_{r\leq}J_{\leq k}^{\equiv i}\,|
> |{}_{r\leq}J_{\leq k}^{\equiv (i+1)}\,|\ \}\,.
$$
If $R$ is empty, then $f_i(J)$ is undefined.  Otherwise, 
$$
f_i(J):=J\setminus r\cup(r\sh+1)\,,\quad\text{where}\quad
r:=\min(R)\,.
$$ 
(ii) Let 
$$
R':=\{\,r'\text{\ \ s.t.\ \ for all $k\leq r'$,}\ \
|{}_{k\leq}J_{\leq r'}^{\equiv (i+1)}\,|
> |{}_{k\leq}J_{\leq r'}^{\equiv i}\,|\ \}\,.
$$
If $R'$ is empty, then $e_i(J)$ is undefined.  Otherwise, 
$$
e_i(J):=J\setminus r'\cup (r'\sh-1)\,,\quad\text{where}\quad
r':=\max(R')\,.
$$ 
}
\\[-1em]
The importance of the crystal operators lies in the following
Refined Demazure Character Formula (cf. Jimbo, et al. \cite{jimbo}).  Define the {\it weight} of a tableau $J\in\CC(L_m)$ by $\wt(L_m):=\Lam_m$
and $\wt(f_i(J)):=\wt(J)-\alpha_i$.
\\[1em]
{\bf Proposition}\quad {\it The character of the Demazure module $V_w(\Lam_m)$ is the weight generating function of the crystal graph $\CC_w(L_m)$:  that is,
$$
%\mathop{\rm char}V_w(\Lam_m):=
\sum_{\mu}\ \dim_\C\! V_w(\Lam_m)_\mu\ \,e^\mu
=\sum_{J\in\CC_w(L_m)} e^{\wt(J)}\,.
$$
In particular,\ $\dim_\C\! V_w(\Lam_m)=\#\CC_w(L_m)$.
}
\\[1em]
Let us give a more pictorial way to understand these operators in the spirit of Lascoux-Schutzenberger \cite{LS}: we progressively remove elements of $J$ which are irrelevant to the action.
We call $j\in J^{\equiv i}$ the {\it $i$-elements} of $J$, and we write sets as usual in increasing order: $J=\{\cdots\sh<j_{-1}\sh<j_{0}\}$.  We start by removing all $j\in J$ except the $i$- and $(i\sh+1)$-elements.  We consider each remaining $i$-element which is immediately followed by an $(i\sh+1)$-element, and we remove these pairs.  
Now we look again for remaining $i$-elements followed by $(i\sh+1)$-elements, and remove these pairs.  
After finitely many iterations, we are left with
a finite subset 
$$J'=\{j'_1<\cdots<j'_s<j''_1<\cdots<j''_t\}\subset J
\quad\text{with all}\ \ j'_k\equiv i\sh+1,\ \
j''_k\equiv i\,.
$$
Then we take $r=j''_1$, the smallest $i$-element, and $r':=j'_s$, the largest $(i\sh+1)$-element of $J'$, so that:
$$
f_i(J)=J\setminus j''_1\cup(j''_1\sh+1)\,,
\qquad
e_i(J)=J\setminus j'_s\cup(j'_s\sh-1)\,.
$$

\noindent {\bf Example}\quad
We exhibit the action of $e_2,f_2$ on the $J$ from our previous example.  This time, we write the elements of $J$ reduced modulo $n=5$:  
since $J$ is $n$-bounded, this loses no information.  
We have underlined the elements to be removed. 
$$
\begin{array}{ccr@{\,}c@{\,}c@{\,}c@{\,}c@{\,}c@{\,}c@{\,}c@{\,}c@{\,}c@{\,}c@{\,}c@{\,}c@{\,}c@{\,}c@{\,}c@{\,}c@{\,}c@{\,}l}
J&=&\{\,\cdots,&\!-3,\,&\!-2,\,&\!-1,\,&\,0,\,&3,\,&4,\,&7,&10,&12,&14,&17,&18,&23,&27,&32,&33,&35,&37\,\}
\\
&=&\cdots\ &2&3&4&5&3&4&2&5&2&4&2&3&3&2&2&3&5&2
\\
&\Rightarrow&\cdots\ &\underline2&\underline3&&&3&&2&&2&&\underline{2}&\underline{3}&3&2&\underline2&\underline3&&2
\\
&\Rightarrow&&&&&&3&&2&&\underline2&&&&\underline3&2&&&&2
\\
J'&=&&&&&&\mathbf3&&\mathbf2&&&&&&&\mathbf2&&&&\mathbf2\\[.5em]
f_2(J)&=&\cdots\ &2&3&4&5&\mathbf3&4&\mathbf3&5&2&4&2&3&3&\mathbf2&2&3&5&\mathbf2 \\
f_2^2(J)&=& \cdots\ &2&3&4&5&\mathbf3&4&\mathbf3&5&2&4&2&3&3&\mathbf3&2&3&5&\mathbf2 \\
f_2^3(J)&=& \cdots\ &2&3&4&5&\mathbf3&4&\mathbf3&5&2&4&2&3&3&\mathbf3&2&3&5&\mathbf3 \\
f_2^4(J)&=&\text{undefined}\hspace{-2.8em}\\[.5em]
e_2(J)&=& \cdots\ &2&3&4&5&\mathbf2&4&\mathbf2&5&2&4&2&3&3&\mathbf2&2&3&5&\mathbf2 \\
e_2^2(J)&=&\text{undefined}\hspace{-2.8em}\\
\end{array}
$$
Note that the irrelevant elements removed from $J$ are the same as those from $e_i(J)$ and $f_i(J)$, so we can easily perform $e_i$ and $f_i$ repeatedly.

In the previous example we computed $\roof(J)=y(L_m)$, 
where:
$$
y=s_{2}\,s_{1}\,s_{3}\,s_{2}\,s_{0}\,
(s_{4}s_{3}\,s_{2}\,s_{1}\,s_{0})^{11}\,s_{4}\,.
$$
By Theorem 1, this means that $J\in\CC_y(L_m)$:
$$
J=f_{2}^{\bdot}\,f_{1}^{\bdot}\,f_{3}^{\bdot}\,f_{2}^{\bdot}\,f_{0}^{\bdot}\,(f_{4}^{\bdot}\,f_{3}^{\bdot}\,f_{2}^{\bdot}\,f_{1}^{\bdot}\,f_{0}^{\bdot})^{11}\,f_{4}^{\bdot}\, L_{14}\ ,
$$
where each $f_i^{\bdot}$ represents some non-negative integer power of $f_i$.  We see from this the comparative rapidity of the roof algorithm in defining and generating Demazure crystals.\quad$\diamond$

\section{Proof of Theorem 1}

For a set $J\in\CC(L_m)$, let $\CC_y(L_m)$ be the unique minimal Demazure crystal containing $J$, and define the {\it ceiling of $J$} to be the extremal element of $\CC_y(L_m)$:
$$
\ceil(J):=y(L_m)\,.
$$
Thus\
$\CC_w(L_m) =\{J\in\CC(L_m)\mid\ceil(J)\leq w(L_m)\}\,,$ and we can restate:
\\[1em]
{\bf Theorem 1}\quad {\it We have\ \ $\roof(J)=\ceil(J)$\ \
for all $J\in\CC(L_m)$.}
\\[1em]
For $J\neq L_m$, we define:
$$
a(J):=\max\{a\mid L_a\subset J\}
\quad\text{and}\quad
r(J):=\min J_{>a(J)}\,.
$$
That is, $a(J)<r(J)$ are the smallest consecutive elements of $J$ which are not consecutive integers.

We let $e^{\max}_i(J)$ denote the result of applying the highest possible power of the raising operator $e_i$ to $J$, and we let:
$$K:=\e_{r-1}^{\max}J\,,$$
where $r:=r(J)$. Observe that $r\sh-1\in K$
(and thus $K\neq J$), since in
$J_{<r}=L_{a(J)}$, the pairs of consecutive entries congruent to $r\sh-1$ and $r$ are irrelevant for the crystal operation.
\\[1em]
\noindent{\bf Ceiling Lemma}\\[.3em] {\it 
(i) For all $J\in\CC(L_m)$, we have: 
$$r(\ceil(J))=r(J)
\quad\text{and}\quad
a(\ceil(J))=a(J)\,.$$
(ii) With $J\neq L_m$ and $K$ as above, we have:
$$
\ceil(J)>\ceil(K)=s_{r-1}\,\ceil(J)\,.
$$
}\\[0em]
{\bf Roof Lemma}\ \ {\it With $J\neq L_m$ and $K$ as above, we have:
$$
\roof(J)>\roof(K)=s_{r-1}\,\roof(J)\,.
$$
}\\[0em]
Assuming these two Lemmas, we can immediately prove Theorem 1 by induction on the quantity:
$$
\height(J):=\sum_{i\leq 0}\,(j_i-i-m)\,,
$$
a sum with finitely many non-zero terms for $J\in\CC(L_m)$.
If $\height(J)=0$, then $J=L_m$ and there is nothing to prove.  Otherwise, $\height(K)<\height(J)$, and we may assume $\roof(K)=\ceil(K)$.  Then the
Roof and Ceiling Lemmas imply:
$$\qquad\qquad
\roof(J)=s_{r-1}\roof(K)=s_{r-1}\ceil(K)=\ceil(J)\,.
\qquad\Box
$$
\\[0em]
{\it Proof of Ceiling Lemma.}\ \ 
We first prove that $a(\ceil(J))\leq a(J)$. 
Let $a:=a(J)$.
Let $\ceil(J)=s_{i_t}\cdots
s_{i_1}$, a reduced decomposition. Then for some
$c_t,\ldots,c_1\geq 0$, $J=f_{i_t}^{c_t}\cdots f_{i_1}^{c_1}L_m$.
The sequence $\{i_t,\ldots,i_1\}$ must contain a subsequence
$\{a\sh+1,\ldots,m\}$. Let $\{j_k,\ldots,j_1\}$ be the rightmost such subsequence: that is, ${j_1}$ is the rightmost occurrence of
$m$ in $\{i_t,\ldots,i_1\}$;  and for $k=1,\ldots,m\sh-a\sh-1$, after ${j_k}$ has been determined
let ${j_{k+1}}$ be the rightmost occurrence of $m\sh-k$ in
$\{i_t,\ldots,i_1\}$ to the left of ${j_k}$.
Let $f_{i}^{\max}T$ denote the result of
applying the lowering operator $f_{i}$ as many times as possible
to $T$; thus, for example, $\ceil(J)=f_{i_t}^{\max}\cdots
f_{i_1}^{\max}L_m$. Then
$$
a(f_{i_t}^{\max}\cdots f_{i_1}^{\max}L_m)\leq 
a(f_{{j_k}}^{\max}\cdots f_{j_1}^{\max}L_m)=m\sh-k,
$$
for $k=1,\ldots,m\sh-a$.  Setting $k=m\sh-a$, we obtain the result.

We prove (i) and (ii) together by induction on $\height(J)$. Let $r:=r(J)$.
If
$\height(J)=0$ then $\ceil(J)=J=L_m$, so (i) is true, and (ii) is
vacuously true.  Assume $\height(J)>0$.  Note that
$$\ceil(J)\geq\ceil(K)\geq s_{r-1}\,\ceil(J)\,,$$ so (ii) is
equivalent to $\ceil(J)\neq\ceil(K)$.

If $r(J)=a(J)+2$, then $a(K)=a(J)+1$.  Since
$\height(K)<\height(J)$, by induction, $a(\ceil(K))=a(K)$, thus
$a(\ceil(K))=a(K)>a(J)\geq a(\ceil(J))$, implying
$\ceil(K)\neq\ceil(J)$.  Therefore $\ceil(J)=s_{r-1}\ceil(K)$, and
clearly (i) follows.

If $r(J)>a(J)+2$, on the other hand, let $wL_m=\ceil(K)$. Since
$\height(K)<\height(J)$, by induction we have $a(wL_m)=a(K)=a$ and 
$r(wL_m)=r(K)=r\sh-1$. Define
$$
v=s_{a+1}s_{a+2}\cdots s_{r-2}\,w\,.
$$ 
Note that $a(vL_m)=a\sh+1$. 
Let $v=s_{i_t}\cdots s_{i_1}$ be a reduced decomposition.  
Then $w=s_{r-2}\cdots s_{a+1}s_{i_t}\cdots s_{i_1}$, 
also a reduced decomposition. Indeed, if we define $k$ by $r(wL_m)=w_k$
(where $wL_m=\{\cdots\sh>w_{-2}\sh>w_{-1}\sh>w_0\}$\,), then $$
(s_{a+(j+1)}\cdots
s_{a+1}s_{i_t}\cdots s_{i_1}L_m)_k = 
1+(s_{a+j}\cdots
s_{a+1}s_{i_t}\cdots s_{i_1}L_m)_k
$$ 
for $j=1,\ldots,r\sh-a\sh-3$, and
$$
(s_{a+1}s_{i_t}\cdots s_{i_1}L_m)_k = 
1+(s_{i_t}\cdots
s_{i_1}L_m)_k\,.
$$
In other words, with each successive
multiplication of $v=s_{i_t}\cdots s_{i_1}$ by $s_{a+j}$ for $j=1,\ldots, r\sh-a\sh-2$, the product increases. 

Now suppose $\ceil(J)=wL_m$.
Then $J=f_{r-2}^{c_{r-2}}\cdots
f_{a+1}^{c_{a+1}}f_{i_t}^{d_t}\cdots f_{i_1}^{d_1}L_m$ for some
$c_{r-2},\ldots,c_{a+1},d_t,\ldots,d_1\geq 0$.  Thus,
$e_{a+1}^{c_{a+1}}\cdots e_{r-2}^{c_{r-2}}J\in C_v(L_m)$.  Thus\\
$a(e_{a+1}^{c_{a+1}}\cdots e_{r-2}^{c_{r-2}}J)\geq a(vL_m)=a\sh+1$, a
contradiction.  Therefore $\ceil(J)\neq wL_m$, so
$\ceil(J)=s_{r-1}wL_m$, from which (i) follows immediately.
\hfill$\Box$
\\[1em]
{\it Proof of Roof Lemma.} For $i\in\mathbb{Z}$, define
$\roof_i(J)$ by
$$\roof_i(J):=\roof(L_i\cup J).$$ Several properties of
$\roof_i(J)$ follow easily from the definition:
\begin{itemize}
\item[1.] $\roof_i(J)=\roof_{i+1}(J)$ \ if $i+1\in J$.
\item[2.] $\roof_i(J)=\roof(\roof_{i+1}(J)\setminus \{i+1\})$ \ if $i+1\not\in J$.
\item[3.] $\roof_i(J)=\roof(J)$ if $J\supset L_i$ ; \ $\roof_i(J)=L_i$ \ if $L_i\supset J$.
\end{itemize}
For $T\in\CC(L_k)$, define
$$\lub(T)=\min_{\geqlex}\left\{J'\in \CC(L_k)\Big|
\begin{array}{l}
J'\geqB T,\\
J'\hbox{ is n-stable}
\end{array}
\right\}.
$$
If $T$ has at most one seam, then $\roof(T)=\lub(T)$. Since
$\roof_{i+1}(J)\setminus (i\sh+1)$ has at most one seam, property
2 above can be modified:
\begin{itemize}
\item[4.] $\roof_i(J)=\lub(\,\roof_{i+1}(J)\setminus (i\sh+1)\,)$ \ if\,
$i\sh+1\not\in J$.
\end{itemize}

For $k\in\mathbb{Z}_{\geq 0}$, let $r[k]:=r+kn$. We will prove the following eight
statements $\bA_k$---$\bH_k$ together by decreasing
induction on $k$. Then we will show that the Roof Lemma is a
consequence of statement $\bC_0$ (i.e., $\bC_k$ for $k=0$).
\begin{itemize}
\item[$(\bA_k)$] Either $\roof_{r[k]}(K)=\roof_{r[k]}(J)$ or
$\roof_{r[k]}(K)=s_{r-1}\roof_{r[k]}(J)$.
\item[$(\bB_k)$] $\roof_{r[k]}(K)\leqB\roof_{r[k]}(J)$.\\[-.5em]
\item[$(\bC_k)$] If $r[k]\in J\setminus K$, $r[k]-1\in K\setminus J$,
then $\roof_{r[k]-2}(K)=s_{r-1}\roof_{r[k]-2}(J)$.\\[-.5em]
\item[$(\bD_k)$] If $r[k]\in J\cap K$, but $r[k]-1\not\in J$ or $K$,
then $\roof_{r[k]-2}(K)=\roof_{r[k]-2}(J)$.\\[-.5em]
\item[$(\bE_k)$] If $r[k], r[k]-1\not\in J\cup K$, then
$\roof_{r[k]-2}(K)=\roof_{r[k]-2}(J)$.\\[-.5em]
\item[$(\bF_k)$] If $r[k]-1\in J\cap K$, but $r[k]\not\in J$ or $K$,
then $\roof_{r[k]-2}(K)=\roof_{r[k]-2}(J)$.\\[-.5em]
\item[$(\bG_k)$] Either $\roof_{r[k]-2}(K)=\roof_{r[k]-2}(J)$ or
$\roof_{r[k]-2}(K)=s_{r-1}\roof_{r[k]-2}(J)$.
\item[$(\bH_k)$] $\roof_{r[k]-2}(K)\leqB\roof_{r[k]-2}(J)$.\\[-.5em]
\end{itemize}
%\\[1em]
%The only other possibility, namely $r[k]\in K\setminus J$,
%$r[k]-1\in J\setminus K$, does not occur, since $K$ is obtained
%from $J$ by applying $e_{r-1}$ several times.
Our induction proof will establish 
the following implications:
\[
(\bA_{k+1}\sh-\bF_{k+1})\ \Rightarrow\ (\bG_{k+1},\bH_{k+1})
\ \Rightarrow\ (\bA_{k},\bB_k)\ \Rightarrow\ (\bC_k,\bD_k,\bE_k)
\]
\[
(\bD_m, \bE_m, \bF_m\!:\,m\sh>k\,)\ \Rightarrow\ (\bF_k)
\]
%\\[.5em]
For the starting point of induction, select $k$ large enough so
that $r[k-2]>j_0$. For such $k$, by property 3,
$\roof_i(K)=\roof_i(L)=L_i$ for $i=r[k\sh-2],r[k\sh-1],r[k]$. Thus
$\bA_k\sh-\bH_k$ are trivially true.
\\[1em]
$(\bA_{k+1}\sh-\bF_{k+1})\ \Rightarrow\ (\bG_{k+1}, \bH_{k+1})$ :  
\\[1em]
Let us restate this as:
$(\bA_{k}\sh-\bF_{k})\Rightarrow(\bG_{k}, \bH_{k})$. If any of the
hypotheses of $\bC_k-\bF_k$ are satisfied, then $\bC_k-\bF_k$
imply $\bG_k$ and $\bH_k$.  There are two
possibilities omitted from the hypotheses of
$\bC_k-\bF_k$:
\begin{itemize}
\item[(i)] $r[k]-1\in J\setminus K$, $r[k]\in K\setminus J$, and
\item[(ii)] $r[k]-1, r[k]\in J\cap K$.
\end{itemize}
However, (i) cannot occur, since $K$ is obtained from $J$ by
applying the raising operator $e_{r-1}$ several times. If (ii)
occurs, then by property 1, $\roof_{r[k-2]}(J)=\roof_{r[k]}(J)$,
and $\roof_{r[k-2]}(K)=\roof_{r[k]}(K)$.  Thus, in this case as
well, $\bG_k$ and $\bH_k$ follow immediately from $\bA_k$ and
$\bB_k$.
\\[1em]
$(\bG_{k+1}, \bH_{k+1})\Rightarrow (\bA_k,\bB_k)$ :
\\[1em]
 Define $t$ by
$\roof_{r[k]}(J)=\up^t(\roof_{r[k+1]-2}(J))$. Then necessarily
$\roof_{r[k]}(K)=\up^t(\roof_{r[k+1]-2}(K))$. Letting
$T=\roof_{r[k+1]-2}(J)$, $U=\roof_{r[k+1]-2}(K)$, we show that
\begin{eqnarray}\label{eqn_TU}
&&\hbox{Either }\ \up^i(U)=\up^i(T)\ \hbox{ or
}\ \up^i(U)=s_{r-1}\up^i(T),\ \hbox{ and}\\
&&\up^i(U)\leqB\up^i(T)
\end{eqnarray}
 $0\leq i\leq t$, by induction on $i$; the
result is then obtained by setting $i=t$.
\\[1em]
We have $\up^{0}(T):=T$, $\up^{0}(U):=U$; thus (1) and (2) hold
for $i=0$. Let $0<i\leq t$, and assume that (1) and (2) hold for
$i-1$. Then either
\\[.5em]
(i) $\up^{i-1}(U)=\up^{i-1}(T)$, in which case
$$\up^i(U)=\up(\up^{i-1}(U))=\up(\up^{i-1}(T))=\up^{i}(T),\hbox{ or }$$
(ii) $\up^{i-1}(U)=s_{r-1}\up^{i-1}(T)$ and $\up^{i-1}(U)\leqB
\up^{i-1}(T)$. In this case, define $p, q$ by
$$\up^i(T)=\up^{i-1}
(T)\setminus p\cup q.$$ Then it is easy to see that
$\up^i(U)=\up^{i-1}(U)\setminus p\cup q'$, where
\[
q'=
\begin{cases}
q, &\hbox{ if }q\not\equiv r-1, r\mod n\\
q+1, &\hbox{ if }q\equiv r-1\mod n\\
q-1, &\hbox{ if }q\equiv r\mod n\\
\end{cases}.
\]
Thus $\up^i(U)=s_{r-1}\up^i(T)$ and $\up^{i}(U)\leqB \up^{i}(T)$.
This proves (1) and (2).
\\[1em]
$(\bD_m, \bE_m, \bF_m: m>k) \Rightarrow(\bF_k)$ :
\\[1em]
Let $m>k$ be the minimum integer such that not both $r[m]-1$ and
$r[m]$ are in $K$. Then, by the definition of the raising operator
$e_{r-1}$, it is not possible that $r[m]\in J\setminus K$,
$r[m]-1\in K\setminus J$. Thus either the hypotheses of $\bD_m$,
$\bE_m$, or $\bF_m$ must hold. Thus
$\roof_{r[m]-2}(J)=\roof_{r[m]-2}(K)$. Since also
$$
_{r[k]-2\leq}J_{\leq r[m]-2}\ =\ _{r[k]-2\leq }K_{\leq r[m]-2}\,,$$
we have
$\roof_{x}(J)=\roof_{x}(K)$, for $r[k]-2\leq x\leq r[m]-2$.
\\[1em]
$(\bA_k,\bB_k)\Rightarrow(\bD_k)$ :
\\[1em]
By Property 1, $\roof_{r[k]-1}(J)=\roof_{r[k]}(J)$ and
$\roof_{r[k]-1}(K)=\roof_{r[k]}(K)$. Now $\bA_k$, $\bB_k$ imply that
if $J'\in\CC(L_m)$ is $n$-stable and $J'\geqB
\roof_{r[k]-1}(K)\setminus \{r[k]-1\}$, then $J'\geqB
\roof_{r[k]-1}(J)\setminus \{r[k]-1\}$ (the same statement with $J$
and $K$ switched holds obviously). The result follows from
Property 4.
\\[1em]
$(\bA_k,\bB_k)\Rightarrow(\bC_k)$ :
\\[1em]
Let
$J'=\mathop{\max}_{\geqB}\{\roof_{r[k]}(J),s_{r-1}\roof_{r[k]}(J)\}$. By
Property 1, $\roof_{r[k]-1}(J)=\roof_{r[k]}(J)$. If $J'\in\T$ is
$n$-stable and $J'\geqB \roof_{r[k]-1}(J)\setminus \{r[k]-1\}$, then
$J'\geqB \hat{J}\setminus \{r[k]-1\}$; conversely, if $J'\in\T$ is
$n$-stable and $J'\geqB \hat{J}\setminus \{r[k]-1\}$, then $J'\geqB
\roof_{r[k]-1}(J)\setminus \{r[k]-1\}$. Thus
$\roof_{r[k]-2}(\hat{J})=\roof_{r[k]-2}(J)$.

We claim that $K\leqB s_{r-1}K$.   Indeed, let $m>k$ be the
minimum integer greater than $k$ such that not both $r[m]-1,r[m]\in
K$. Then by the definition of the raising operator $e_{r-1}$, it
is not possible that $r[m]\in K$. Thus $\roof_{r[m]-2}(K)\leqB
s_{r-1}\roof_{r[m]-2}(K)$.  Define $t'$ by
$\roof_{r[k]}(K)=\up^{t'}(\roof_{r[m]-2}(K))$.  Then it is easy to
see that $\up^{i}(\roof_{r[m]-2}(K))\leqB
s_{r-1}\up^{i}(\roof_{r[m]-2}(K))$, $1\leq i\leq t'$. This proves
the claim.

Define $t$ by $\roof_{r[k]-2}(\hat{J})=\up^t(\hat{J})$.  Then by
Property 1,
$\roof_{r[k]-2}(K)=\roof_{r[k]-1}(K)=\up^t(\roof_{r[k]}(K))$. For
$1\leq i\leq t$, let
$$\up^i(\hat{J})=\up^{i-1}
(\hat{J})\setminus p\cup q.$$ Then
$\up^i(\roof_{r[k]}(K))=\up^{i-1}(\roof_{r[k]}(K))\setminus p\cup
q'$, where
\[
q'=
\begin{cases}
q, &\hbox{ if }q\not\equiv r \mod n\\
q-1, &\hbox{ if }q\equiv r\mod n\\
\end{cases}.
\] The result follows from this.
\\[1em]
$(\bA_k,\bB_k)\Rightarrow(\bE_k)$ :
\\[1em]
We have that
$\roof_{r[k]-2}(J)=\roof(\roof_{r[k]}(J)\setminus\{r[k]-1,r[k]\})$,
$\roof_{r[k]-2}(K)=\roof(\roof_{r[k]}(K)\setminus\{r[k]-1,r[k]\})$. Let
$L=\roof_{r[k]}(J)\setminus\{r[k]-1,r[k]\}$.  If $L$ has only one seam, then
the result is obvious. Thus assume
that $L$ has two seams:
\begin{eqnarray*}
S_1&=&\{r[k]=p_1<\cdots<p_t\}\\
S_2&=&\{r[k]-1=p_{t+1}<\cdots<p_{t+s}\}
\end{eqnarray*}
where $p_{i+1}=p_i+n$, $1\leq i\leq t+s-1$, $i\neq t$. Then
$\up^t(L)$ has exactly one seam, namely $S_2$ with possibly some
additional elements added to its tight end. Thus
$\roof_{r[k]-2}(J)=\roof(\up^t(L))=\lub(\up^t(L))$. We claim that
$\lub(\up^t(L))=\lub(L)$. The claim implies
$\roof_{r[k]-2}(J)=\lub(L)=\lub(\roof_{r[k]}(J)\setminus\{r[k]-1,r[k]\})$;
replacing $J$ with $K$, in precisely the same manner we show that
$\roof_{r[k]-2}(K)=\lub(\roof_{r[k]}(K)\setminus\{r[k]-1,r[k]\})$. But
it is clear that $\lub(\roof_{r[k]}(J)\setminus\{r[k]-1,r[k]\})=
\lub(\roof_{r[k]}(K)\setminus\{r[k]-1,r[k]\})$.  Thus the result
follows from the claim.

To prove the claim, note that since $\up^t(L)\geqB L$,
$\lub(\up^t(L))\geqlex\lub(L)$. It suffices to show that
$\lub(L)\geqB \up^t(L)$, since this implies
$\lub(L)\geqlex\lub(\up^t(L))$. We show something slightly
stronger: if $M\geqB L$ is $n$-stable, then $M\geqB \up^t(L)$. We
can express $M=L\setminus\{p_1,\ldots,p_{t+s}\}\cup
\{p_1,\ldots,p_{t+s}\}$, where $x_i> p_i$, $i=1,\ldots,t+s$.
Likewise, $\up^t(L)=L\setminus\{p_1,\ldots,p_{t}\}\cup
\{q_1,\ldots,q_{t}\}$, where $\up^i(L)=\up^{i-1}(L)\setminus
p_i\cup q_i$. To show $M\geqB \up^t(L)$, it suffices to show that
$q_i\leq x_i$, $i=1,\ldots,t$. This is clear from the definition
of the $\up$ operation. Indeed, let $i_{\max}$ be the largest $i$
for which $q_i-p_i> n$. There are no tight ends in $\up^t(L)$
between $r[k]$ and $q_{i_{\max}}$; thus $q_i\leq x_i$, $1\leq i\leq
i_{\max}$. If $i_{\max}<t$, then for $i_{\max}<i\leq t$,
$$q_i=\min\{q\not\in L\mid q-p_i\leq n-1,\,
q-n\in\up^{i-1}(L),\, q\not\equiv r\mod n,\}.$$ Inductively, this
implies that $q_i\leq x_i$.
\\[2em]
This completes the proof of $\bA_k-\bH_k$. Noting that $r[0]=r$, we see that
$\bC_0$ implies $\roof_{r-2}(K)=s_{r-1}\roof_{r-2}(J)$.  By
property 3, $\roof(J)=\roof_{a(J)}(J)$ and
$\roof(K)=\roof_{a(J)}(K)$. Using identical arguments as in the
proof of $(\bG_{k+1}, \bH_{k+1})\Rightarrow (\bA_k,\bB_k)$, we see
that $\roof_{a(J)}(K)=s_{r-1}\roof_{a(J)}(J)$, which completes the
proof of the Roof Lemma. \hfill$\Box$

\section{Proof of Proposition 3}

\subsection{Proof of Proposition 3(i)}
The result states that:
$$
v_J=\sum_K a_K^J\, \ep_K\,,
$$
where the sum runs over $K\leqlex J$.
We use induction on $\ell$, where $\roof(J)=\up^\ell(J)$.  If $\ell=0$, then $v_J=\ep_J$ and there is nothing to prove.  

Now let $\ell>0$. We inductively apply the Proposition to
$
\hatJ :=\up(J)=J\setminus p\cup q\,,
$\ \ so that:
$$
v_\hatJ = \sum_\hatK a_\hatK^\hatJ \,\ep_\hatK\,,
$$
where the sum runs over $\hatK\leqlex 
\hatJ$.  Thus:
$$
v_J=\hatE_{pq}\, v_\hatJ 
=\sum_\hatK \sum_{h\in\Z}
a_\hatK^\hatJ\, E_{p+nh,q+nh}\,\ep_{\hatK}\,.
$$
It suffices to show the following:
\\[1em]
{\bf Lemma}\quad {\it
Let $\hatJ:=\up(J)=J\setminus p\cup q$. Consider any $\hatK\leqlex\hatJ$ and any $h\in\Z$ such that:
$$
p':=p\sh+nh\not\in\hatK\quad\text{and}\quad 
q':=q\sh+nh\in\hatK\,.
$$
Then we have:
$$
K:=\hatK\setminus q'\cup p'\ \leqlex\ J=\hatJ\setminus q\cup p\,.
$$
}
\\
We prove the Lemma using several facts which follow easily from the definitions.  Let $K,J\in\CC(L_m)$.
For $J=\{\cdots\sh<j_{-1}\sh<j_0\}$ with $\ord(J)=m$, recall that:\ \
$
\height(J):=\sum_{i\leq 0}\,(j_i-i-m)\,.
$
\begin{enumerate}

\item  If $a_K^J\neq 0$, then $\height(K)=\height(J)$.

\item  If $a_K^J\neq 0$,\ 
then $|K^{\equiv i}_{>N}| = |J^{\equiv i}_{>N}|$ for all $N\ll0$.

\item  If $\hatJ=\up(J)=J\setminus p\cup q$,\
then $\hatJ^{\equiv q}\subset \hatJ_{\leq q}$.

\item  
%If $\hatJ=\up(J)=J\setminus p\cup q$
If $J_{>p}$ contains no loose ends of $J$ (that is, $j\sh-n\in J$ for all $j\in J_{>p}$),
\ \ and\ \
$|K^{\equiv i}_{\geq p}|=
|J^{\equiv i}_{\geq p}|$ for some $i$, then
$K^{\equiv i}_{\geq p} \geqB
J^{\equiv i}_{\geq p}$.
\end{enumerate}

Proceeding with the proof of the Lemma, suppose first that $\hatK=\hatJ$.  By Fact 3 we must have $q'\leq q$, so $p'\leq p$ and clearly
$K=\hatJ\setminus q'\cup p'\,\leqlex\,
\hatJ\setminus q\cup p =J$.

Now let $\hatK\llex\hatJ$, and let $\hatk$ be the {\it split point}, the value such that:
$$
\hatk\in\hatK,\ \ \hatk\not\in \hatJ, 
\quad\text{and}\quad
K_{<\hatk}=J_{<\hatk}.
$$

Case (a):  $\hatk<p$.  Then: $$
\hatk\in\hatK,\ \ \hatk\not\in J
\quad\text{and}\quad \hatK_{<\hatk}=\hatJ_{<\hatk}=J_{<\hatk}\,,
$$
so $\hatK\leqlex J$.  But clearly $K\leqlex \hatK$, so $K\leqlex J$ as desired.

Case (b):  $p\leq\hatk$. If $p'\leq p$, then clearly
$K\leqlex J$ as desired.  
On the other hand, suppose $p<p'$.  
Then $K_{<p}=\hatK_{<p}=\hatJ_{<p}=J_{<p}$ by the definition of $\hatk$.  
Furthermore, for all $i$ and some $N<p$ we have
$|K^{\equiv i}_{>N}|=|J^{\equiv i}_{>N}|$ by Fact 2, 
and thus 
$|K^{\equiv i}_{\geq p}|=|J^{\equiv i}_{\geq p}|$. By definition $J_{>p}$ contains no loose ends, so Fact 4 implies that $K^{\equiv i}_{\geq p}\geqB J^{\equiv i}_{\geq p}$ for all $i$.  We also have $K_{<p}=J_{<p}$, so $K\geqB J$.  If $K\stackrel{\rm B}{>}J$, then clearly $\height(K)>\height(J)$, contradicting Fact 1.  We conclude that
$K=J$, and we are done.  

This proves the Lemma, and hence Proposition 3(i).

\subsection{Proof of Proposition 3(ii)}
To derive the formula relating the leading coefficients $a_J^J$ and $a_\tJ^\tJ$, note first that
$$\begin{array}{rcl}
v_J&=&(\hatE_{p_1q_1}\cdots\hatE_{p_tq_t})v_\tJ\\[.5em]
&=&\displaystyle\sum_K\sum_{h_1,\ldots,h_t\in\Z\!\!\!\!\!}\!\!\!\!
a_K^\tJ\, (E_{p_1+nh_1,q_1+nh_1}\cdots E_{p_t+nh_t,q_t+nh_t})\,\ep_K\\[1.5em]
&=&\displaystyle\sum_K\sum_{h_1,\ldots,h_t\in\Z\!\!\!\!\!}\!\!
\pm\,a_K^\tJ\ \ep_{K{\uparrow}(h_1,\ldots,h_t)}\ ,
\end{array}$$ 
where we use notation:
$$
K{\uparrow}(h_1,\ldots,h_t):=K\setminus\{q_1\sh+nh_1
,\ldots ,q_t\sh+nh_t\}\cup\{p_1\sh+nh_1,\cdots
p_t\sh+nh_t\}
$$
provided $q_i\sh+nh_i\in K$ and $p_i\sh+nh_i\not\in K$ for all $i\leq t$ ; otherwise $K{\uparrow}(h_1,\ldots,h_t)$ is undefined, and $\ep_{K{\uparrow}(h_1,\ldots,h_t)}:=0$.
\\[1em]
{\bf Lemma} {\it (i) If $J=K{\uparrow}(h_1,\ldots,h_t)$ for some $K$ with $a_K^\tJ\neq 0$, then $K=\tJ$.
\\
(ii) If $J=\tJ{\uparrow}(h_1,\ldots,h_t)$, then there is a unique permutation $\sigma$ of $\{1,2,\ldots,r\}$ such that
$$
p_i\sh+nh_i=p_{\sigma(i)}\,,\qquad q_i\sh+nh_i=q_{\sigma(i)}\,.
$$
We obtain in this way every permutation $\sigma$ satisfying
$q_i-p_i=q_{\sigma(i)}-p_{\sigma(i)}$ for all $i\leq t$.
}
\\[1em]
The Proposition follows easily from (i) and (ii) of the Lemma, since:
$$\begin{array}{rcl}
v_J&=&\displaystyle\sum_K\!\!\!\!
\mathop{\sum_{h_1,\ldots,h_t\in\Z\!\!\!\!\!}}
_{K{\uparrow}(h_1,\ldots,h_t)=J}\!\!\!\!
a_K^\tJ\, (E_{p_1+nh_1,q_1+nh_1}\cdots E_{p_t+nh_t,q_t+nh_t})\,\ep_K
\\[3em]
&\stackrel{(i)}=&\left(a_\tJ^\tJ\!\!\!\!
\displaystyle
\mathop{\sum_{h_1,\ldots,h_t\in\Z\!\!\!\!\!}}
_{\tJ{\uparrow}(h_1,\ldots,h_t)=J}\!\!\!\!
(E_{p_1+nh_1,q_1+nh_1}\cdots E_{p_t+nh_t,q_t+nh_t})
\,\ep_{\tJ}\right)+
\text{lower}
\\[3em]
&\stackrel{(ii)}=&\left(a_\tJ^\tJ\,
\displaystyle\sum_{\sigma}
(E_{p_{\sigma(1)},q_{\sigma(1)}}\cdots E_{p_{\sigma(r)},q_{\sigma(r)}})
\,\ep_{\tJ}\right)+
\text{lower}\\[2em]
&\stackrel{(*)}=&\left(a_\tJ^\tJ\,
\displaystyle
\mathop{\sum_{\sigma}}
(E_{p_1q_1}\cdots E_{p_tq_t})
\,\ep_{\tJ}\right)+
\text{lower}\\[2em]
&=&\pm\, a_\tJ^\tJ\cdot\#\{\sigma\}\cdot\ep_J\ +\ \text{lower}\,,
\end{array}$$
where $\sigma$ runs over the set of all permutations of $\{1,\ldots,t\}$ such that $q_i-p_i=q_{\sigma(i)}-p_{\sigma(i)}$ for $i\leq t$ : clearly $\#\{\sigma\}=\prod_{d\geq 1} \mu_d!$ .
Equation $(*)$ holds because the operators $E_{p_iq_i}$ all commute for $i\leq t$.  It remains to prove the Lemma.
\\[1em]
{\it Proof of Lemma (i).}  Suppose
$J=K{\uparrow}(h_1,\ldots,h_t)=\tJ{\uparrow}(0,\ldots,0)$ and $a_K^{\tJ}\neq 0$. Let
$$
p':=\min\{p_1\sh+nh_1,\ldots,p_t\sh+nh_t\}\,,\qquad
p:=p_1=\min\{p_1,\ldots,p_t\}\,.
$$
We clearly have $K_{<\min(p,p')}=\tJ_{<\min(p,p')}$.
If $p'<p$, then $p'\not\in K$, $p'\in\tJ$, and $K_{<p'}=\tJ_{<p'}$, so $\tJ\llex K$, which contradicts Proposition 3(i).

Thus $p\leq p'$, and $K_{<p}=\tJ_{<p}$.  Furthermore, by Fact 2 in the proof of Prop. 3(i), for any $i\leq t$ we have $|K^{\equiv i}_{>N}|=|\tJ^{\equiv i}_{>N}|$ for some $N<p$.  Hence for any $i$,\ \ $|K^{\equiv i}_{\geq p}|=|\tJ^{\equiv i}_{\geq p}|$ .  Since $\tJ_{>p}$ clearly has no loose ends, Fact 4 implies  $K^{\equiv i}_{\geq p}\geqB\tJ^{\equiv i}_{\geq p}$ for any $i$, and we also know $K_{<p}=\tJ_{<p}$.  We conclude that $K\geqB\tJ$, and a fortiori $K\geqlex\tJ$.  Since $K\leqlex\tJ$ by Proposition 3(i), we must have $K=\tJ$.
\\[1em]
{\it Proof of Lemma (ii).}  Suppose
$J=\tJ{\uparrow}(0,\ldots,0)=\tJ{\uparrow}(h_1,\ldots,h_t)$.  Define:
$$
p'_i:=p_i+nh_i\,,\qquad q'_i:=q_i+nh_i\,,\qquad
d_i:=q_i-p_i:=q'_i-p'_i$$
Then we have
$\{p_1,\ldots,p_t\}=\{p'_1,\ldots,p'_t\}$
and $\{q_1,\ldots,q_t\}=\{q'_1,\ldots,q'_t\}$,
so there exist permutations $\alpha,\beta$ of $\{1,\ldots,r\}$ such that:
$$
p_i = p'_{\alpha(i)}\,,\qquad
q_i = q'_{\beta(i)}\,.
$$
We will use the following facts:
\begin{enumerate}
\item We have $p_i=p+n(i-1)$ for $i=1,\ldots,r$,\, and also\ $q_1<q_2<\cdots<q_t$ .

This follows from the seam-pulling action of the up-operation.
\item If $i<j$ and $q_i\equiv q_j\mod n$, then
$d_i\geq d_j$ .

Indeed, for a fixed $k$, the set of all 
$q_i\equiv k\mod n$ forms an arithmetic progression 
$\{q,q\sh+n,\ldots\}$, whereas the corresponding set of $\{p_i\mid q_i\equiv k\}$ is a subset of the arithmetic progression $\{p_1, p_2,\ldots\}=\{p,p\sh+n,\ldots\}$.  Hence, if $i<j$ and $q_i\equiv q_j\mod n$, then $p_j-p_i\geq q_j-q_i$, and so $d_i=q_i-p_i\geq q_j-p_j=d_j\,.$

\item  $\alpha(i)=\beta(i)
\quad\Longleftrightarrow\quad
d_i=d_{\alpha(i)}
\quad\Longrightarrow\quad
q_i\equiv q_{\alpha(i)}\mod n\quad
$

If $\alpha(i)=\beta(i)$, then
$
d_i=q_i-p_i=q'_{\beta(i)}-p'_{\alpha(i)}
=q'_{\alpha(i)}-p'_{\alpha(i)}=d_{\alpha(i)}
\,,$
and also $q_i=q'_{\beta(i)}\equiv q_{\beta(i)}=q_{\alpha(i)}$.

%\item $ q_i\equiv q_\beta(i)\equiv q_{\beta^{-1}(i)}\mod n$.
%\quad
%This is because $q_i=q'_{\beta(i)}\equiv q_{\beta(i)}$.
\end{enumerate}
\smallskip

\noindent
Assume $\alpha\neq\beta$, and let $j$ be the smallest value such that $\alpha(j)\neq \beta(j)\,.$\
Then $j$ is minimal with $\beta^{-1}\alpha(j)\neq j$, and necessarily:
$$
\beta^{-1}\alpha(j)>j\quad\text{with}\quad
q_{\beta^{-1}\alpha(j)}=q'_{\alpha(j)}\equiv q_{\alpha(j)}\,.
$$
 
Consider the sequence:\ 
$
j,\, \alpha(j),\, \alpha^2(j),\, \alpha^3(j),\ldots\,.
$\ \
If $\alpha(j),\alpha^2(j),\cdots,\alpha^{c}(j)\sh<j$, then by the definition of $j$ 
%we have $\alpha^2(j)=\beta\alpha(j)$,
%$\alpha^3(j)=\beta\alpha^2(j)=\beta^2\alpha(j)$, etc., 
and Fact 3 we have:
$$\begin{array}{c}
d_j\neq d_{\alpha(j)}=d_{\alpha\alpha(j)}=d_{\alpha\alpha\alpha(j)}=\cdots=d_{\alpha^{c+1}(j)}\\[.7em]
q_{\alpha(j)}\equiv q_{\alpha\alpha(j)}\equiv
q_{\alpha\alpha\alpha(j)}\equiv\cdots\equiv q_{\alpha^{c+1}(j)}\,.
\end{array}$$
But we eventually have $\alpha^{c+1}(j)=j$, so to avoid the contradiction $d_j\neq d_j$, there must exist some $k:=\alpha^{c+1}(j)$ such that:
$$
k>j\quad\text{with}\quad
d_k=d_{\alpha(j)}
\quad\text{and}\quad
q_k\equiv q_{\alpha(j)}\,.
$$

Case (a): $j<\beta^{-1}\alpha(j) \leq k$.
Then by Fact 2, we have $d_{\beta^{-1}\alpha(j)}\geq d_k=d_{\alpha(j)}$.  But:
$$\begin{array}{rcl}
d_{\beta^{-1}\alpha(j)}&=&
q_{\beta^{-1}\alpha(j)}-p_{\beta^{-1}\alpha(j)}\\
&<&
q_{\beta^{-1}\alpha(j)}-p_{j}\\
&<&
q'_{\alpha(j)}-p'_{\alpha(j)}=d_{\alpha(j)}
\,,\end{array}$$
so this case is impossible.

Case (b): $j<k<\beta^{-1}\alpha(j)$.
Then by Fact 1, we have $p_j<p_k<q_k<q_{\beta^{-1}\alpha(j)}$.
But:
$$\begin{array}{rcl}
q_{\beta^{-1}\alpha(j)}
&=&q'_{\alpha(j)}\\
&=&p'_{\alpha(j)}+d_{\alpha(j)}\\
&=&p_j+d_{k}\\
&<&p_k+d_k\ =\ q_k\,.
\end{array}$$
Thus, this case is impossible also.

The above contradictions show that $\alpha=\beta$.  Hence we have
$$
p+nh_i=p'_i=p_{\sigma(i)}\,,\qquad
q+nh_i=q'_i=q_{\sigma(i)}\,,
$$
where $\sigma=\alpha^{-1}=\beta^{-1}$, which is the first part of Lemma (ii).

To see the second part of Lemma (ii), suppose $p_i,q_i$ given and let $\sigma$ satisfy
$d_i=d_{\sigma(i)}$.  Then define $h_i:=(p_{\sigma(i)}-p_i)/n$, so that
$p'_i:=p_i+nh_i=p_{\sigma(i)}$ and:
$$
q'_i:=q_i+nh_i=p'_i+d_i=p_{\sigma(i)}+d_{\sigma(i)}=q_{\sigma(i)}\,.
$$
Thus $\{p_1,\ldots,p_t\}=\{p'_1,\ldots,p'_t\}$
and $\{q_1,\ldots,q_t\}=\{q'_1,\ldots,q'_t\}$,
so $J=\tJ{\uparrow}(h_1,\ldots,h_t)$, as desired.

This proves the Lemma, and hence Proposition 3(ii).

\end{document}